\documentclass[a4paper,11pt,reqno]{amsart}
\usepackage{geometry} 
\usepackage{mathtools}
\usepackage{graphicx}
\usepackage{paralist}
\usepackage[utf8]{inputenc}
\usepackage[T1]{fontenc}
\usepackage{csquotes}
\usepackage[all]{xy}
\usepackage{chngcntr}
\usepackage{amsthm}
\usepackage{amsmath}

\usepackage{amsfonts}
\usepackage{amssymb}
\usepackage{mathrsfs}
\usepackage{MnSymbol}
\usepackage{tikz-cd}
\usepackage{tikz, tikz-cd}
\usepackage{mathdots}
\usepackage{stmaryrd}
\usepackage{hyperref}
\usepackage{fancyhdr}
\usepackage{xcolor}

\newtheorem{theorem}{Theorem}
\newtheorem{thm}{Theorem}[section]
\newtheorem{lem}[thm]{Lemma}
\newtheorem{prop}[thm]{Proposition}
\newtheorem{cor}[thm]{Corollary}
\newtheorem{rem}[thm]{Remark}
\theoremstyle{remark}

\theoremstyle{claim}

\newtheorem*{thm*}{Theorem}

\theoremstyle{definition}

\theoremstyle{example}

\theoremstyle{convention}

\counterwithin{equation}{section}

\theoremstyle{convention}

\newcommand{\bbN}{\mathbb{N}}

\newcommand{\bbR}{\mathbb{R}}

\newcommand{\bbH}{\mathbb{H}}

\DeclareMathOperator{\Ai}{Ai}
\DeclareMathOperator{\Bi}{Bi}
\newcommand{\etalchar}[1]{$^{#1}$}

\author[J.~Clutterbuck]{Julie~Clutterbuck}
\address{School of Mathematics \\
Monash University \\
9 Rainforest Walk,
VIC 3800 Australia}
\email{Julie.Clutterbuck@monash.edu}

\author[F.~J\"ackel]{Frieder~J\"ackel}
\address{Département de Mathématiques \\ Universit\'e Libre de Bruxelles \\ Campus de la Plaine, CP 210, Boulevard du Triomphe, B-1050 Bruxelles, Belgique}
\email{frieder.jackel@ulb.be}

\author[X.~H.~Nguyen]{Xuan~Hien~Nguyen}
\address{Department of Mathematics \\ Iowa State University \\ 411 Morrill Rd, Ames, IA 50011 USA}
\email{xhnguyen@iastate.edu \textrm{(Corresponding author)}}

\date{\today}
\thanks{The research of Julie Clutterbuck was supported by grants DP220100067 and DP250103808 of the Australian Research Council. Frieder~J\"ackel was supported by the DFG priority program "Geometry at infinity" and the grant FNRS EoS O.0033.22F of the Fonds de la Recherche Scientifique (F.R.S.-FNRS)}

\calclayout

\begin{document}
\title[Constant potentials do not minimise the fundamental gap in $\mathbb{H}^n$]{Constant potentials do not minimise the fundamental gap on convex domains in hyperbolic space}

\begin{abstract}
We show that for every $n \geq 2$ and $D > 0$ there exist a convex domain $\Omega \subseteq \bbH^n$ with diameter $D$ and a convex potential $V$ on $\Omega$ such that the fundamental gap of the operator $-\Delta+V$ is strictly smaller than the fundamental gap of $-\Delta$. In comparison to previous work, this result requires more refined control of the eigenfunctions.
\end{abstract}

\maketitle

\tableofcontents

%

\section{Introduction}

We consider the low Dirichlet eigenvalues of the Laplace operator $- \Delta + V(x)$ with convex potential $V$ on a convex, compact domain $\Omega \subseteq \mathbb H^n$. This operator has a discrete spectrum with $\infty$ as its accumulation point. If the sequence of eigenvalues is arranged in increasing order $\lambda_1 < \lambda_2 \leq \lambda_3 \leq \cdots$, the \emph{fundamental gap} $\Gamma(\Omega;V)$ is the difference between the first two eigenvalues,
\[
     \Gamma(\Omega;V) \coloneq \lambda_2 - \lambda_1>0.
\]
If the potential is constant, we simply write $\Gamma(\Omega)$. The main result of the present article is the following.

\begin{theorem}\label{Main Theorem}
For every $n \geq 2$ and every $D > 0$ there exist a convex domain $\Omega \subseteq \bbH^n$ and a convex function $V$ on $\Omega$ such that ${\rm diam}(\Omega)=D$ and $\Gamma(\Omega;V) < \Gamma(\Omega)$.
\end{theorem}

\subsubsection*{Motivation}
A flurry of activity has focused on the first part of the Fundamental Gap Conjecture, i.e., whether $\Gamma(\Omega)$ or $\Gamma(\Omega;V)$ can be bounded from below in terms of the diameter. The present article is the first one to study the second part of the Fundamental Gap Conjecture, i.e., if $\Gamma(\Omega;V)$ is, among all convex potentials, minimised by the constant ones on manifolds other than $\bbR^n$. The only investigation of this question after \cite{fundamental} was in the case of graphs.

For one-dimensional space, Lavine \cite{Lavine-gap} proved that constant potentials maximise the fundamental gap among all convex potentials. Recent work by El Allali and Harrell \cite{AH22} found that the constant potential is not optimal among single-well potentials. 

For graphs, Ahrami et al \cite{ahrami2024optimizing} showed that on a metric tree graph the optimal potential (in the class of convex potentials) must be piecewise affine.  However they uncovered an interesting phenomenon:    there are graphs where the constant potential does not minimise. They suspect this property is generic. Metric graphs display some aspects of the one-dimensional case -- locally, we are dealing with ordinary differential equations -- and some aspects of the manifold case, in that the domain can have non-trivial topology, and super-linear growth of length enclosed in a geodesic ball can be understood as a curvature condition.

\subsection{Background}\label{subsec:background}
In their paper on the Fundamental Gap Conjecture, Andrews and the first author  \cite{fundamental} showed that for every bounded convex domain $\Omega \subseteq \bbR^n$ and every convex potential $V$ on $\Omega$ we have
\[
    \frac{3\pi^2}{{\rm diam}(\Omega)^2} \leq \Gamma(\Omega) \leq \Gamma(\Omega;V).
\]
The result is sharp, with the limiting case being rectangles that collapse to a line. We refer to the introduction of \cite{fundamental} and the survey by Dai--Seto--Wei \cite{daifundamental} for more information about the history and earlier work on this subject.

The Fundamental Gap Conjecture in other spaces of constant sectional curvature is more difficult to investigate.  Dai, He, Seto, Wang, and Wei (in various subsets) \cite{SWW19,DSW21,he2017fundamental} generalized  the  estimate to convex domains in $\mathbb{S}^n$, showing that the same bound $\Gamma(\Omega;V) \geq 3\pi^2/{\rm diam}(\Omega)^2$ holds (though the condition for non-constant potentials is stronger than just convexity - see \cite[Equation (1.4)]{CWY25}). We also note the probabilistic proof of  Cho--Wei--Yang \cite{CWY25}.
 
Together with Bourni--Stancu--Wei--Wheeler the first and third authors \cite{BCN+21, BCN+22} showed that the bound for Euclidean and spherical spaces does \emph{not} hold in hyperbolic space. Namely, they showed that for all $n \geq 2$ and all positive constants $D,\varepsilon>0$ there exists a convex domain $\Omega \subseteq \bbH^n$ with diameter $D$ whose fundamental gap satisfies $\Gamma(\Omega) < \varepsilon \frac{ \pi^2}{D^2}$. This was then generalized by Khan and the third author \cite{KN24} to all negatively curved Hadamard manifolds. However, it was shown by Khan--Saha--Tuerkoen \cite{KST25} and Khan--Tuerkoen \cite{KT24} that for \emph{horo-convex} domains $\Omega \subseteq \bbH^n$ one still has $\Gamma(\Omega) \geq c(n,D)>0$. An alternative proof for horo-convex domains with small diameter was given by Wei--Xiao \cite{WX25}.

\subsection{Strategy}\label{subsec: strategy}

The central idea here is variational:   if a potential $V_0$ minimises the fundamental gap among all convex potentials, then for a family of  convex potentials $V_t$  
\[ \left. \frac{d}{dt} \Gamma(\Omega,V_t)\right|_{t=0}=0.\]
Here we will rely on standard results on the analyticity of eigenvalues under perturbation of the operator.   Moreover, we have the very concrete Hellmann--Feynman formula for this derivative, 
\[ \left. \frac{d}{dt} \Gamma(\Omega,V_t)\right|_{t=0}= \int_\Omega \left.\frac{d V_t }{dt}\right|_{t=0} (u_2^2-u_1^2) \, d{\rm vol}\]
where $u_1,u_2$ are the first and second eigenfunctions associated to $-\Delta + V_0$. 

If we can find a potential $V_t=V_0+ tP$, where $V_0$ is constant (without loss of generality, $V_0\equiv 0$) and $P$ is convex, such that 
\begin{equation}\label{eq: core of strategy}
	\int_\Omega P (u_2^2-u_1^2) \, d{\rm vol} < 0,
\end{equation}
then we can conclude that for small $t > 0$, the fundamental gap for the potential $V_t=tP$ is strictly smaller than that for the zero potential. 

We implement the argument,  firstly by choosing a convex region $\Omega$ which, in the upper half-space model for $\mathbb{H}^2$, is amenable to separation of variables.   The separation of variables allows us to rather explicitly construct the eigenfunctions as solutions to coupled ordinary differential equations.   We choose a simple distance function as the potential $P$. 

The main new insight of this paper is that  we show that \eqref{eq: core of strategy} holds  through fine estimates on the solutions to the ODEs for  the eigenfunctions $u_1$ and $u_2$.    After appropriate rescalings,  the relevant ordinary differential operator turns out to be a perturbation of the Airy operator $-\frac{d^2}{dx^2}+x$. Therefore, by exploiting perturbation arguments, we will show that  the eigenfunctions behave essentially  as the eigenfunctions of the Airy operator on $[0,\infty)$. From this we  deduce the Theorem \ref{Main Theorem}.

Many of these arguments will be familiar from earlier work.   The variational argument was used to opposite effect by Lavine \cite{Lavine-gap}, who also used comparisons with the Airy operator.  The separation of variables construction can be found in Shih \cite{shih1989counterexample}.  The domain is similar to those of \cite{BCN+21, BCN+22}, but inspired by the  metric graphs example of Ahrami et al \cite[Example 5.9]{ahrami2024optimizing}, we use only half the domain.  The metric graphs example also motivates the choice of $P$.   

However the  control on the eigenfunctions that we require is significantly more refined than in the previous work. 
Namely, the proof of \cite[Theorem 1.1]{BCN+22} essentially only needs that the first eigenfunction is exponentially small in the so-called neck region. This is certainly not enough to control the integral in (\ref{eq: core of strategy}), which is why we 
use the approximation by the Airy eigenfunctions.

\subsection{Structure of the article.}
The article is structured as follows. In Section \ref{subsec: domain and potential} we define the domains and the potential for which we will verify (\ref{eq: core of strategy}), and explain the separation of variables ansatz that reduces the PDE eigenvalue problem to an ODE eigenvalue problem. Section \ref{subsec: Airy equation} recalls the necessary background material on the Airy equation. In Section \ref{subsec: perturbation estimates} we collect the required elementary perturbation estimates. The heart of the article is contained in Section \ref{subsec: Airy approx}, in which we show that, after appropriate rescalings, the ODE eigenvalue problem from Section \ref{subsec: domain and potential} is a small perturbation of the Airy eigenvalue problem, and use this to approximate the relevant eigenfunctions by the eigenfunctions of the Airy equation. The proof of Theorem \ref{Main Theorem} in the case $n=2$ is then presented in Section \ref{subsec: proof of main thm}. Finally, in Section \ref{sec: high dim}, we explain how to adapt the arguments to general $n \geq 2$.

\vspace{1cm}

\noindent
\subsection*{Acknowledgements.} This research originated at Simons-Laufer Mathematical Sciences Institute in Fall 2024, during the program \emph{New Frontiers in Curvature: Flows, General Relativity, Minimal Submanifolds, and Symmetry}. We thank the Institute for giving us the opportunity to meet and work on the problem. 
The second author also wants to thank Bernhard Kepka for helpful conversations about perturbation theory. The third author thanks Mustard Seed Community Farm for their hospitality during part of this work.


\section{Set up}\label{sec: Set up}

\subsection{Choice of domain and potential}\label{subsec: domain and potential}
In this subsection we introduce the precise convex domains $\Omega$ and convex potential $P$ for which we will verify (\ref{eq: core of strategy}). For the sake of clarity, we work in dimension 2 here then discuss in Section \ref{sec: high dim} how to generalize to higher dimensions. 

We work with the upper half-space model  $\bbH^2 = \{(x,y) \in \bbR^2 \, | \, y > 0\}$ of hyperbolic space. For all $\varphi_0 \in (0,\frac{\pi}{2})$ and $\mu > 0$ we consider the convex region
\begin{equation}\label{eq: def of Omega}
	\Omega \coloneq \Omega_{\varphi_0,\mu} \coloneq \Big\{\big(r\sin(\varphi),r\cos(\varphi)\big) 
    | \, r \in (1,e^{\pi/\sqrt{\mu}}) \text{ and } \varphi \in (0,\varphi_0) \Big\}.
\end{equation}
Observe that in these polar coordinates for $\bbH^2$, $\{\varphi=0\} 
$ is a 
geodesic. Therefore,
\begin{equation}\label{eq: def of potential}
	P \colon \Omega 
    \to \bbR, \, p \mapsto d_{\bbH^2}\big(p,\{\varphi=0\}\big)
\end{equation}
is convex as $\bbH^2$ is non-positively curved (see for example \cite[Theorem 1.3 on p.~4]{BGS85}). Moreover, because dilations are isometries in the upper half-space model, in the polar coordinates $P$ only depends on $\varphi$ so  we can write $P(\varphi)$. In fact, the only properties of $P$ that we will use are that $P$ depends only on $\varphi$ and $P^\prime(\varphi) > 0$ everywhere.

The reason for the choice (\ref{eq: def of Omega}) of $\Omega$ is that it allows for a separation of variables ansatz, and thus for a reduction to an ODE (see Lemma \ref{lem: strategy} below). We review the relevant parts for the convenience of the reader, and refer to \cite[Section 2]{BCN+21} for further details.

In the polar coordinates from (\ref{eq: def of Omega}), the hyperbolic Laplacian is given by
\[
	\Delta u = \cos^2(\varphi)\partial_{\varphi \varphi}^2(u)+r^2\cos^2(\varphi)\partial_{rr}^2(u)+r\cos^2(\varphi)\partial_{r}(u).
\]
We write
\[
	u(r,\varphi)=f(r)h(\varphi).
\]
Then, solving $-\Delta u = \lambda u$ is equivalent to simultaneously solving
\[
	r^2f^{\prime \prime}(r)+rf^\prime(r)=-\tilde{\mu}f
	\quad \text{and} \quad
	h^{\prime \prime}(\varphi)=\big(\tilde{\mu}-\lambda \cos^{-2}(\varphi)\big)h(\varphi),
\]
where $\tilde{\mu} \in \bbR$ is an auxiliary parameter. The equation for $f$ has a solution on $(1,e^{\pi/\sqrt{\mu}})$ with Dirichlet boundary conditions if and only if $\tilde{\mu}=j^2 \mu$ for some $j \in \bbN$, in which case the solution is (up to scaling) given by $f_j(r)=\sin(j \sqrt{\mu} \log(r))$. 

The equation for $h$ is now
\begin{equation*}
	h^{\prime \prime}(\varphi)=\big(j^2 \mu- \lambda\cos^{-2}(\varphi)\big)h (\varphi) \text{ for } \varphi \in (0,\varphi_0)
	\quad \text{and} \quad
	h(0)=h(\varphi_0)=0,
\end{equation*} 
which for each $j\in\mathbb{N}$, has solutions $h=h^{(j^2\mu)}_k$ associated with $\lambda=\lambda^{(j^2\mu)}_k$, where $k\in \mathbb{N}$.  The doubly-indexed collection $\lbrace \lambda^{(j^2\mu)}_k\rbrace_{j,k\in \mathbb{N}}$  are the eigenvalues of $-\Delta$ on $\Omega$.  We are only concerned with the two smallest.    The smallest is found when $j=k=1$, so $\lambda_1=\lambda^{(\mu)}_1$.  For the second smallest eigenvalue of $-\Delta$, we will claim that $\lambda_{2}^{(\mu)}<\lambda_{1}^{(4\mu)}$, so that  $\lambda_2=\lambda_{2}^{(\mu)}$, and only the $j=1$ case is relevant for our investigation.    The  eigenfunctions are given by $u_k(r,\varphi)=f_1(r)h_{k}^{(\mu)}(\varphi),$ for $k=1,2$. See \cite[Section 2.3]{BCN+21} for more details.

Therefore, one is left with analyzing the eigenvalue problem for $h$:
\begin{equation}\label{eq: eq for h}
	h^{\prime \prime}(\varphi)=\big(\mu-\lambda \cos^{-2}(\varphi)\big)h(\varphi) \text{ for } \varphi \in (0,\varphi_0)
	\quad \text{and} \quad
	h(0)=h(\varphi_0)=0.
\end{equation}

Note that the volume form is given by $d{\rm vol}_{\rm hyp}=r^{-1}\cos^{-2}(\varphi) dr \wedge d\varphi$ in the polar coordinates from (\ref{eq: def of Omega}). Therefore, the following lemma summarizes the above discussion. In its formulation, we use the notations from above.

\begin{lem}[Strategy]\label{lem: strategy}
Let $\varphi_0 \in (0,\frac{\pi}{2})$ and $\mu > 0$. Denote by $h_k$ the $k$-th eigenfunction of the eigenvalue problem 
\begin{equation*}
	h^{\prime \prime}(\varphi)=\big(\mu-\lambda \cos^{-2}(\varphi)\big)h(\varphi) \text{ for } \varphi \in (0,\varphi_0)
	\quad \text{and} \quad
	h(0)=h(\varphi_0)=0
\end{equation*}
with the weighted $L^2$-normalization
\begin{equation}\label{eq: L^2 normalization h}
	\int_0^{\varphi_0} \cos^{-2}(\varphi)h_k^2(\varphi) \, d\varphi=1.
\end{equation}
Assume
\begin{equation}\label{eq: strategy eigenvalue assumption}
	\lambda_2^{(\mu)} < \lambda_1^{(4\mu)}
\end{equation}
and
\begin{equation}\label{eq: core of strategy 2}
	\int_0^{\varphi_0} P(\varphi)\left(h_2^2(\varphi)-h_1^2(\varphi) \right)\cos^{-2}(\varphi)  \, d\varphi < 0,
\end{equation}
where $P$ denotes the potential defined in (\ref{eq: def of potential}). Then, for $t > 0$ small enough, the fundamental gap of the convex potential $V_t=tP$ in $\Omega_{\varphi_0,\mu}$ satisfies 
\[
    \Gamma(\Omega_{\varphi_0,\mu} ; V_t) < \Gamma(\Omega_{\varphi_0,\mu}),
\]
where $\Omega_{\varphi_0,\mu} \subseteq \bbH^2$ is the convex domain defined in (\ref{eq: def of Omega}).
\end{lem}

We will show that, for $\varphi_0 \in (0,\frac{\pi}{2})$ fixed, the assumptions (\ref{eq: strategy eigenvalue assumption}),(\ref{eq: core of strategy 2}) will be satisfied if $\mu$ is chosen large enough. In fact, we will obtain much better estimates. Indeed, (\ref{eq: eigenvalue asymptotic expansion}) and Corollary \ref{cor: integral expansion} give precise asymptotics for the eigenvalues and the value of the integral as $\mu \to \infty$. We have, for definiteness, chosen the specific $P$ from (\ref{eq: def of potential}). However, Corollary \ref{cor: integral expansion} will show that (\ref{eq: core of strategy 2}) will hold if $\mu$ is large enough for every potential $P$ with $P^\prime(\varphi_0)>0$. Moreover, the proof shows that other eigenvalue gaps can be decreased also (see Remark \ref{rem: other gaps}).

To achieve this, we require a much better understanding of the behaviour of the eigenfunctions $h_k$ compared to \cite{BCN+22}. Namely, the proof of \cite[Theorem 1.1]{BCN+22}  only needs that $h_1$ (and its derivative) are exponentially small in the "neck region" $\{\varphi \approx 0\}$ as $\mu \to \infty$. This is certainly not enough to control the integral in (\ref{eq: core of strategy 2}). We will instead show that, as $\mu \to \infty$, the eigenfunctions $h_k$ are essentially appropriate rescalings of the eigenfunctions of the Airy equation on $[0,\infty)$ with Dirichlet boundary conditions.

\subsection{The Airy equation}\label{subsec: Airy equation}

We use this subsection to list some facts about the Airy equation, introduce notation, and evaluate the integral involving the eigenfunctions of the Airy equation that will serve as the asymptotic model for the integral in (\ref{eq: core of strategy 2}).

The \emph{Airy equation} is
\[
	y^{\prime \prime}(x)=xy(x).
\]
The Airy functions $\Ai$ and $\Bi$ are distinct and linearly independent solutions of the Airy equation. Both are defined for all $x \in \bbR$, both are oscillating for $x < 0$, and for $x > 0$ $\Ai$ is exponentially decaying while $\Bi$ is exponentially growing (see for example \cite[Figure 3.1 on p.~69, and (3.5.17a),(3.5.17b) on p.~100]{BO99}). We denote the zeros of $\Ai$ by
\[
	\dots < -a_3 < -a_2 < -a_1 < 0.
\]

If $y:[0,\infty) \to \bbR$ is a non-zero solution of the Airy eigenvalue problem
\[
	y^{\prime \prime}(x)=(x-\alpha) y(x)
\]
with $y(0)=0$ and $y \in L^2(0,\infty)$, then $\Ai(-\alpha)=0$, i.e., $\alpha=a_k$ for some $k \in \bbN$. Indeed, any solution $y$ of $y^{\prime \prime}(x)=(x-\alpha) y(x)$ is a linear combination of $\Ai(\cdot-\alpha)$ and $\Bi(\cdot-\alpha)$. Since $y \in L^2(0,\infty)$ and $\Bi(x)$ grows exponentially for $x \to \infty$, $y$ is a multiple of $\Ai(\cdot - \alpha)$; whence $\Ai(-\alpha)=0$ follows from $y(0)=0$.

Therefore, for all $k \in \bbN$, 
\begin{equation}\label{eq: def of Airy eigfct on half-line}
	v_k(x) \coloneq \frac{\Ai(x-a_k)}{||\Ai(\cdot - a_k)||_{L^2(0,\infty)}}
\end{equation}
is the $L^2$-normalized $k$-th eigenfunction of the Airy operator $-\frac{d^2}{dx^2}+x$ on $[0,\infty)$ with Dirichlet boundary conditions and eigenvalue $a_k$.

In Corollary \ref{cor: integral expansion}, we will use the following as  the (negative of the) asymptotic model for the desired integral in (\ref{eq: core of strategy 2}).

\begin{lem}[Model Integral]\label{lem: model integral}
We have
\[
	\int_0^\infty x\left(v_2^2(x)-v_1^2(x) \right) \, dx = \frac{2}{3}(a_2-a_1) > 0,
\]
where $v_1,v_2$ are defined in (\ref{eq: def of Airy eigfct on half-line}), and $0 > -a_1 > -a_2 > \dots$ are the zeros of $\Ai$.
\end{lem}

\begin{proof}
For each $\gamma > 0$, take the differential operator $A_\gamma \coloneq -\frac{d^2}{dx^2}+\gamma x$, and denote by $\alpha_k(\gamma)$ ($k \in \bbN$) the $k$-th eigenvalue of $A_\gamma$ considered as an unbounded operator on $L^2(0,\infty)$ with Dirichlet boundary conditions (see the proof of Lemma \ref{lem: Airy on finite and infinite domain} for further functional analytic details).

Observe that for any $L^2$-normalized eigenfunction $v$ of the Airy operator with eigenvalue $\alpha$, the function $v_\gamma(x)=\gamma^{1/6}v(\gamma^{1/3}x)$ is an $L^2$-normalized eigenfunction of $A_\gamma$ with eigenvalue $\gamma^{2/3}\alpha$. Thus, $\alpha_k(\gamma)=\gamma^{2/3}a_k$ for all $k \in \bbN$. In particular, $\alpha_k^\prime(1)=\frac{2}{3}a_k$. Together with the Hellmann--Feynman identity  we obtain
\[
	\frac{2}{3}a_k=\alpha_k^\prime(1)=\big(\left(\partial_\gamma A_\gamma|_{\gamma=1}\right)v_k,v_k \big)_{L^2(0,\infty)}=\int_0^\infty x v_k^2(x) \, dx.
\]
This implies the desired identity.
\end{proof}

\section{Approximation by the Airy equation}\label{sec: Approximation}

In this section, we will show that the eigenvalue problem (\ref{eq: eq for h}) is, after appropriate rescaling, a perturbation of the Airy eigenvalue problem. This will allow us to control the eigenvalues and the integral (\ref{eq: core of strategy 2}) in terms of the Airy equation. 

To achieve this, we recall some elementary perturbation estimates that bound the difference of the eigenvalues and eigenvectors of two operators in terms of the difference of the two operators. This is carried out in Section \ref{subsec: perturbation estimates}. Then, in Section \ref{subsec: Airy approx}, we introduce the appropriate rescalings of the eigenvalue problem (\ref{eq: eq for h}) allowing for an application of the results from Section \ref{subsec: perturbation estimates}. Finally, we present the proof of Theorem \ref{Main Theorem} in Section \ref{subsec: proof of main thm}.

\subsection{Preliminary perturbation estimates}\label{subsec: perturbation estimates}

Let $\mathcal{H}=(H,\langle \cdot , \cdot \rangle_{\mathcal{H}})$ and $\widetilde{\mathcal{H}}=(H,\langle \cdot , \cdot \rangle_{\widetilde{\mathcal{H}}})$ be two real separable Hilbert space structures on a common underlying vector space $H$ whose norms are equivalent, i.e., there exists $C_0 \geq 1$ such that
\begin{equation} \label{eq: equivalent norms}
    \frac{1}{C_{0}}||v||_{\mathcal{H}} 
    \leq ||v||_{\widetilde{\mathcal{H}}}
    \leq C_0 |||v||_{\mathcal{H}} 
    \quad \text{for all }v \in H.
\end{equation}
Let $A$, resp.~$\tilde{A}$, be an unbounded self-adjoint operator on $\mathcal{H}$, resp.~$\widetilde{\mathcal{H}}$, with compact resolvent. Denote by $(u_i)_{i \in \bbN} \subseteq \mathcal{H}$, resp.~$(\tilde{u}_i)_{i \in \bbN} \subseteq \widetilde{\mathcal{H}}$, a complete orthonormal basis consisting of eigenvectors of $A$, resp.~$\tilde{A}$, with eigenvalues $\alpha_1 \leq \alpha_2 \leq \dots$, resp.~$\tilde{\alpha}_1 \leq \tilde{\alpha}_2 \leq \dots$. We also assume that $u_i \in {\rm Dom}(\tilde{A})$ and $\tilde{u}_i \in {\rm Dom}(A)$ for all $i \in \bbN$.

\begin{lem}\label{lem: perturbation}
There exists a constant $C$ only depending on $C_0$ with the following properties. Let $k \in \bbN$ and assume $\varepsilon_k,\tilde{\varepsilon}_k \in \bbR_{\geq 0}$ are such that
\begin{equation}\label{eq: perturbation norm distortion 1}
    (1+\varepsilon_k)^{-1}||u||_{\mathcal{H}}^2
    \leq ||u||_{\widetilde{\mathcal{H}}}^2
    \leq (1+\varepsilon_k)||u||_{\mathcal{H}}^2
    \quad \text{for all } u \in {\rm span}\{u_1,\dots,u_k\}
\end{equation}
and
\begin{equation}\label{eq: perturbation norm distortion 2}
    (1+\tilde{\varepsilon}_k)^{-1}||\tilde{u}||_{\widetilde{\mathcal{H}}}^2
    \leq ||\tilde{u}||_{\mathcal{H}}^2
    \leq (1+\tilde{\varepsilon}_k)||\tilde{u}||_{\widetilde{\mathcal{H}}}^2
    \quad \text{for all } \tilde{u} \in {\rm span}\{\tilde{u}_1,\dots,\tilde{u}_k\}.
\end{equation}
Then we have
\begin{equation}\label{eq: eigenvalue perturbation}
	- C\left(\sum_{i=1}^k \big|\big|(\tilde{A}-A)\tilde{u}_i \big|\big|_{\mathcal{H}}^2 \right)^{\frac{1}{2}}
    -\tilde{\varepsilon}_k \tilde{\alpha}_k
    \leq \tilde{\alpha}_k - \alpha_k \leq \varepsilon_k \alpha_k+ C\left(\sum_{i=1}^k \big|\big|(\tilde{A}-A)u_i \big|\big|_{\mathcal{H}}^2 \right)^{\frac{1}{2}}.
\end{equation}
Moreover, if the $k^{th}$ eigenvalue $\alpha_k$ of $A$ is simple, then, after potentially changing $u_k$ up to sign, we have
\begin{equation}\label{eq: eigenvector perturbation}
	||\tilde{u}_k -u_k|| \leq \frac{C}{\Gamma_k}
	\left(\varepsilon_k \alpha_k+\left(\sum_{i=1}^k \big|\big|(\tilde{A}-A)u_i \big|\big|_{\mathcal{H}}^2 \right)^{\frac{1}{2}}+\tilde{\varepsilon}_k\tilde{\alpha}_k+
	\left(\sum_{i=1}^k \big|\big|(\tilde{A}-A)\tilde{u}_i \big|\big|_{\mathcal{H}}^2 \right)^{\frac{1}{2}}
	\right),
\end{equation}
where $\Gamma_k \coloneq \min\{\alpha_{k+1}-\alpha_k,\alpha_{k}-\alpha_{k-1}\} > 0$.
\end{lem}

We believe that this is well-known, however we were unable to locate a precise reference, and so we provide the proof for the convenience of the reader.

\begin{proof}
We abbreviate $E \coloneq \tilde{A}-A$ and $S_k \coloneq {\rm span}\{u_1, \dots, u_k\}$. Then, by the min-max principle, the Cauchy--Schwarz inequality and the triangle inequality, we get
\begin{align*}
	\tilde{\alpha}_k 
	\leq 
	\max_{u \in S_k} \frac{\langle \tilde{A}u, u \rangle_{\widetilde{\mathcal{H}}}}{||u||_{\widetilde{\mathcal{H}}}^2} 
	&\leq 
	\max_{u \in S_k} \frac{||Au||_{\widetilde{\mathcal{H}}}}{||u||_{\widetilde{\mathcal{H}}}}
	+
	\max_{u \in S_k} \frac{||Eu||_{\widetilde{\mathcal{H}}}}{||u||_{\widetilde{\mathcal{H}}}} \\
	&\leq 
    (1+\varepsilon_k)\max_{u \in S_k} \frac{||Au||_{\mathcal{H}}}{||u||_{\mathcal{H}}}
	+ C_0^2\max_{u \in S_k} \frac{||Eu||_{\mathcal{H}}}{||u||_{\mathcal{H}}} \\
    &=
    (1+\varepsilon_k)\alpha_k
	+ C_0^2\max_{u \in S_k} \frac{||Eu||_{\mathcal{H}}}{||u||_{\mathcal{H}}},
\end{align*}
where in the third inequality we used \eqref{eq: equivalent norms}, \eqref{eq: perturbation norm distortion 1} and $A(S_k) \subseteq S_k$. 
Moreover, we see that $||Eu||_{\mathcal{H}} \leq ||u||_{\mathcal{H}}(\sum_{i=1}^k||Eu_i||_{\mathcal{H}}^2)^{\frac{1}{2}}$ when applying the Cauchy-Schwarz inequality in $\bbR^k$. 

This yields the upper bound in (\ref{eq: eigenvalue perturbation}). Exchanging the roles of $\alpha_k$ and $\tilde{\alpha}_k$, and keeping in mind that $||\cdot||_{\mathcal{H}}$ and $||\cdot||_{\widetilde{{\mathcal{H}}}}$ are equivalent, we also obtain the lower bound.

It remains to prove (\ref{eq: eigenvector perturbation}). Decompose $\tilde{u}_k=\beta u_k +w \in {\rm span}\{u_k\} \oplus {\rm span}\{u_k\}^\perp$. After potentially changing $u_k$ up to sign, we may without loss of generality assume that $\beta\geq 0$. Then $||\tilde{u}_k-u_k|| \leq \sqrt{2}||w||$, and thus it suffices to bound $||w||$.

Note $A\tilde{u}_k=\tilde{A}\tilde{u}_k-E\tilde{u}_k=\tilde{\alpha}_k\tilde{u}_k-E\tilde{u}_k$, whence 
\[
	(A-\alpha_k)w=(A-\alpha_k)\tilde{u}_k=(\tilde{\alpha}_k-\alpha_k)\tilde{u}_k-E\tilde{u}_k.
\]
Now $||w||_{\mathcal{H}} \leq \frac{1}{\Gamma_k}||(A-\alpha_k)w||_{\mathcal{H}}$ since by assumption the restriction of $A-\alpha_k$ to ${\rm span}\{u_k\}^\perp$ has spectral gap $\Gamma_k$. Therefore, (\ref{eq: eigenvector perturbation}) follows by the triangle inequality and (\ref{eq: eigenvalue perturbation}).
\end{proof}

We apply Lemma \ref{lem: perturbation} to compare the eigenvalues and eigenfunctions of $A_R$, the Airy operator on large but finite intervals $(0,R)$ to the eigenvalues and eigenfunctions of $A_{\infty}$, the Airy operator in Section \ref{subsec: Airy equation} on $(0,\infty)$. We will show the difference is essentially irrelevant for $R$ large.

For $R > 0$, we denote by $\alpha_1^R \leq \alpha_2^R \leq \dots$ the eigenvalues and by $(u_k^R)_{k \in \bbN}$ the orthonormal basis of eigenfunctions of the Airy equation on $(0,R)$ with Dirichlet boundary condition, that is,
\[
	\left(-\frac{d^2}{dx^2}+x \right)u_k^R(x)=\alpha_k^R u_k^R(x) \text{ for } x \in (0,R)
	\quad \text{and} \quad
	u_k^R(0)=u_k^R(R)=0
\]
with 
\[
	(u_i^R,u_j^R)_{L^2(0,R)}=\delta_{ij}
	\quad \text{for all } i,j \in \bbN.
\]
Recall from (\ref{eq: def of Airy eigfct on half-line}) that $v_k$ is the $L^2$-normalized eigenfunction of the Airy operator with eigenvalue $a_k$ on $(0,\infty)$ with Dirichlet boundary condition.

\begin{lem}\label{lem: Airy on finite and infinite domain}
For all $k \in \bbN$ we have, after potentially changing $u_k^R$ up to sign,
\[
	\lim_{R \to \infty}\alpha_k^R=a_k
	\quad \text{and} \quad
	\lim_{R \to \infty}\left|\left|u_k^R-v_k|_{[0,R]}\right|\right|_{L^2([0,R])}=0.
\]
Moreover, the convergence is exponentially fast.
\end{lem}

This is well-known, but we provide a proof for the convenience of the reader because we were not able to locate a precise reference in the literature. The reason why Lemma \ref{lem: Airy on finite and infinite domain} is true is because $v_k(R)$ is exponentially small as $R \to \infty$.

\begin{proof}
The Airy operator
\[
	A_\infty \colon \, {\rm Dom}(A_\infty) 
    \to L^2(0,\infty), \, v \mapsto -v^{\prime \prime}+x v,
\]
where 
\[
	{\rm Dom}(A_\infty) \coloneq \left\{v \in L^2(0,\infty) \cap H_{\rm loc}^2(0,\infty)\, | \, -v^{\prime \prime}+xv \in L^2(0,\infty) \text{ and } v(0)=0 \right\},
\] is unbounded, and via \cite[Theorem 9.6]{Te09},
is self-adjoint. Moreover, \cite[Theorem XIII.67]{RS78} shows that $A_\infty$ has purely discrete spectrum with a complete set of eigenfunctions. From the discussion in Section \ref{subsec: Airy equation}, the eigenfunctions of $A_\infty$ are exactly the $v_k$ defined in (\ref{eq: def of Airy eigfct on half-line}) with eigenvalues given by zeroes of the Airy equation $a_k$. 

Similarly,
\[
	A_R \colon \, {\rm Dom}(A_R) 
    \to L^2(0,R), \,
	u \mapsto -u^{\prime \prime}+xu,
\]
where ${\rm Dom}(A_R) \coloneq H^2(0,R) \cap H_0^1(0,R)$, is an unbounded self-adjoint operator with purely discrete spectrum and a complete set of eigenfunctions.

Fix a smooth cut-off function $\rho: \bbR \to [0,1]$ such that $\rho(x)=1$ if $x \leq 0$, $\rho(x)=0$ if $x \geq 1$, and $|\rho^\prime|, |\rho^{\prime \prime}| \leq 10$. Now, for $R \geq 10$ and $i \in \bbN$, we define $v_i^R \in {\rm Dom}(A_R)$ by
\[
	v_i^R(x) \coloneq \rho(x-R+1 )v_i( x).
\]
Then, as $v_i$ and its derivatives are exponentially decaying, we have
\[
	||(A_R-a_i)v_i^R||_{L^2(0,R)} \leq \varepsilon_R
	\quad \text{and} \quad
	(v_i^R,v_j^R)_{L^2(0,R)}=\delta_{ij}+O(\varepsilon_R),
\]
where $\varepsilon_R \to 0$ exponentially fast as $R \to \infty$. Thus, due to the min-max principle, we can bound
\[
	\alpha_k^R \leq \max_{u \in {\rm span} \{v_1^R, \dots, v_k^R\}}\frac{(A_R u,u)_{L^2(0,R)}}{||u||_{L^2(0,R)}^2} \leq a_k+O(\varepsilon_R).
\]

One can show using the Sturm comparison theorem that, for all $R > 0$ and $i \in \bbN$, the eigenfunction $u_i^R$ is exponentially decaying for $x$ large enough (only depending on an upper bound for $\alpha_i^R$); see Claim 1 in the proof of Proposition \ref{prop: main approx} for more details.   Hence, the functions $u_i^\infty \in {\rm Dom}(A_\infty)$ defined by
\[
	u_i^\infty(x) \coloneq u_i^{R,\infty}(x) \coloneq \rho(x-R+1)u_i^R(x)
\]
also satisfy
\[
	||(A_\infty-\alpha_i^R)u_i^\infty||_{L^2(0,\infty)} \leq \varepsilon_R
	\quad \text{and} \quad
	(u_i^\infty,u_j^\infty)_{L^2(0,\infty)}=\delta_{ij}+O(\varepsilon_R),
\]
for some (potentially different) $\varepsilon_R$ with $\varepsilon_R \to 0$ exponentially fast as $R \to \infty$. Again by the min-max principle,
\[
	a_k \leq \max_{v \in {\rm span} \{u_1^\infty, \dots, u_k^\infty\}}\frac{(A_\infty v,v)_{L^2(0,\infty)}}{||v||_{L^2(0,\infty)}^2} \leq \alpha_k^R+O(\varepsilon_R).
\]
Therefore, $|\alpha_k^R-a_k|=O(\varepsilon_R)$. In particular, for $R$ large enough, all eigenvalues of $A_R$ are simple and the distances between consecutive eigenvalues are bounded from below uniformly in $R$. 

The arguments from the proof of Lemma \ref{lem: perturbation} now also show $||v_k^R-u_k^R||_{L^2(0,R)}=O(\varepsilon_R)$. This finishes the proof since $||v_k|_{[0,R]}-v_k^R||_{L^2(0,R)}$ is also exponentially small.
\end{proof}

\subsection{The eigenvalue problem as a perturbed Airy equation}\label{subsec: Airy approx}

We now consider appropriate rescalings of the eigenvalue problem (\ref{eq: eq for h}) in order to interpret it as a perturbed Airy eigenvalue problem; this will allow an application of Lemma \ref{lem: perturbation}. We remark that the change of variables in \eqref{eq: def of tilde(u)} is $\varphi=\varphi_0 - \delta^{1/3} x$ therefore, while the interesting behaviour for $\tilde u_k$ is near $x=0$, the one for $h_k$ is near $\varphi=\varphi_0$. 

Fix $\varphi_0 \in (0,\frac{\pi}{2})$ and take $\mu > 0$ large. Recall that $\lambda_1 < \lambda_2 \leq \dots$ are the eigenvalues and $h_k$ ($k \in \bbN$) the eigenfunctions of problem 
\begin{equation}\tag{\ref{eq: eq for h}}
	h_k^{\prime \prime}(\varphi)=\big(\mu-\lambda_k \cos^{-2}(\varphi)\big)h_k(\varphi) \text{ for } \varphi \in (0,\varphi_0)
	\quad \text{and} \quad
	h_k(0)=h_k(\varphi_0)=0,
\end{equation}
with the weighted $L^2$-normalization
\begin{equation}\tag{\ref{eq: L^2 normalization h}}
	\int_0^{\varphi_0} \cos^{-2}(\varphi)h_k^2(\varphi) \, d\varphi=1.
\end{equation}
We are interested in the behaviour of $\lambda_k$ and $h_k$ as $\mu \to \infty$. To that purpose, we define
\begin{equation}\label{eq: def of delta}
	\delta \coloneq \frac{\mu^{-1}}{2\tan(\varphi_0)},
\end{equation}
and
\begin{equation}\label{eq: def of tilde(u)}
	\tilde{u}_k \colon [0,\varphi_0\delta^{-1/3}] \longrightarrow \bbR, \, x \longmapsto n_k h_k(\varphi_0-\delta^{1/3}x),
\end{equation}
where $n_k \in \bbR_{>0}$ is the unique renormalization constant such that
\begin{equation}\label{eq: L2 normalization tilde(u)}
    \int_0^{\varphi_0 \delta^{-1/3}} \frac{\cos^2(\varphi_0)}{\cos^2(\varphi_0- \delta^{1/3}x)} \tilde{u}_k(x)^2 \, dx=1.
\end{equation}
Explicitly, $n_k$ is given by
\begin{equation}\label{eq: norm approx}
    \frac{\delta^{1/3}}{n_k^2}=\cos^2(\varphi_0).
\end{equation}

\newcommand{\coeffutilde}{\tilde{c}_k} 
\newcommand{\coeffu}{{c}_k} 
\newcommand{\coeffh}{d_k} 

The functions $\tilde u_k$ satisfy the eigenvalue problem for an operator $\tilde A$ close to the Airy operator. We now compute $\tilde A$. Using the notation 
$h''_k(\varphi) =   \coeffh(\varphi) h_k(\varphi)$
for \eqref{eq: eq for h} and 
$\tilde u''_k(x) = \coeffutilde(x) \tilde u_k(x)$
\begin{align}\label{eq: coeff of rescaled eigenfct}
  \coeffutilde(x) 
     & =  \delta^{2/3}\coeffh(\varphi_0-\delta^{1/3}x) \qquad\text{via \eqref{eq: def of tilde(u)}} \notag\\
    &= \delta^{2/3} \left(\mu-\frac{\lambda_k}{\cos^2(\varphi_0-\delta^{1/3}x)} \right) \qquad\text{via \eqref{eq: eq for h}} \notag\\
    &= \delta^{2/3}\mu \left(1-\frac{\cos^2(\varphi_0)}{\cos^2(\varphi_0-\delta^{1/3}x)}-\frac{\cos^2(\varphi_0)}{\cos^2(\varphi_0-\delta^{1/3}x)}\frac{\mu^{-1}\lambda_k-\cos^2(\varphi_0)}{\cos^2(\varphi_0)} \right) \notag\\
    &= \frac{\delta^{-1/3}}{2\tan(\varphi_0)}\left(1-\frac{\cos^2(\varphi_0)}{\cos^2(\varphi_0-\delta^{1/3}x)} \right)-\frac{\cos^2(\varphi_0)}{\cos^2(\varphi_0-\delta^{1/3}x)}\tilde{\alpha}_k,
\end{align}
where we use \eqref{eq: def of delta} in the last line, and define $\tilde{\alpha}_k$ by 
\begin{equation}\label{eq: def of tilde(lambda)}
	\tilde{\alpha}_k \coloneq \frac{\delta^{-1/3}}{2\tan(\varphi_0)}\frac{\mu^{-1}\lambda_k-\cos^2(\varphi_0)}{\cos^2(\varphi_0)}.
\end{equation}
Finally, setting  
\begin{equation}\label{eq: def of perturbation operator}
	\tilde{A} \coloneq \frac{\cos^2(\varphi_0-\delta^{1/3}x)}{\cos^2(\varphi_0)}\left(-\frac{d^2}{dx^2}+\frac{\delta^{-1/3}}{2\tan(\varphi_0)}\left(1-\frac{\cos^2(\varphi_0)}{\cos^2(\varphi_0-\delta^{1/3}x)} \right) \right),
\end{equation}
we have that $\tilde u_k$ satisfies the eigenvalue problem
\(
	\tilde{A}\tilde{u}_k=\tilde{\alpha}_k \tilde{u}_k
\)
for all $k \in \bbN$. 

Observe that $\tilde{A}$ is an unbounded self-adjoint operator on the weighted $L^2$-space 
\[
    \tilde{L}^2\big(0,\varphi_0 \delta^{-1/3}\big) \coloneq \left(L^2\big(0,\varphi_0\delta^{-1/3}\big), \langle \cdot , \cdot \rangle_{\tilde{L}^2} \right)
\]
with the weighted $L^2$-inner product $\langle \cdot , \cdot \rangle_{\tilde{L}^2}$ given by 
\begin{equation}\label{eq: weighted L^2-norm}
    \langle f , g \rangle_{\tilde{L}^2} \coloneq \int_0^{\varphi_0 \delta^{-1/3}} \frac{\cos^2(\varphi_0)}{\cos^2(\varphi_0 - \delta^{1/3}x)} f(x)g(x) \, dx.
\end{equation}
Therefore, due to the weighted normalization (\ref{eq: L2 normalization tilde(u)}), the $\tilde{u}_k$ ($k \in \bbN$) form a complete orthonormal basis for $\tilde{L}^2\big(0,\varphi_0 \delta^{-1/3}\big)$.

The following is the main result of this subsection. It compares the eigenvalues and the eigenfunctions of $\tilde A$ with the corresponding quantities and functions for the Airy operator. For its formulation, we recall that $0>-a_1 > -a_2 > \dots$ are the zeros of $\Ai$, and that $v_k$ denotes the eigenfunction on $(0,\infty)$ of the Airy operator with eigenvalue $a_k$ defined in (\ref{eq: def of Airy eigfct on half-line}).

\begin{prop}\label{prop: main approx}
For every $\varphi_0 \in (0,\frac{\pi}{2})$ and $k \in \bbN$ there exist constants $C> 0$ and $\mu_0 > 0$ such that for every $\mu \geq \mu_0$ we have the following estimates:
\begin{enumerate}[(i)]
	\item (Eigenvalue)
		\begin{equation}\label{eq: eigenvalue approx}
		|\tilde{\alpha}_k-a_k| \leq C \delta^{1/3}
		\end{equation}
	\item (Weighted integral)
		\begin{equation}\label{eq: int approx by Airy}
			\int_0^{\varphi_0 \delta^{-1/3}}x \left|\tilde{u}_k^2(x)-v_k^2(x) \right| \, dx 				\leq C \delta^{1/3}
		\end{equation}
\end{enumerate}
\end{prop}

Note that (\ref{eq: eigenvalue approx}) translates into the following asymptotic expansion for the eigenvalue
\begin{equation}\label{eq: eigenvalue asymptotic expansion}
	\frac{\lambda_k}{\mu}=\cos^2(\varphi_0)\left(1+2\tan(\varphi_0)a_k  \delta^{1/3}+O\big(\delta^{2/3}\big)\right)
	\quad \text{as } \mu \to \infty.
\end{equation}
In particular, the first assumption (\ref{eq: strategy eigenvalue assumption}) in Lemma \ref{lem: strategy} is satisfied for all $\mu$ large enough. The fact that the second assumption (\ref{eq: core of strategy 2}) in Lemma \ref{lem: strategy} is satisfied for all $\mu$ large enough follows from Corollary \ref{cor: integral expansion} below.

The following technical remark will be useful in the proof of Theorem \ref{Main Theorem}.

\begin{rem}\label{rem: locally uniform constants}\normalfont
Given $k \in \mathbb N$, the constants $C$ and $\mu_0$ in Proposition \ref{prop: main approx} can be chosen to be locally uniform in $\varphi_0$, i.e., for all $0 < \varphi_{-} \leq \varphi_{+} < \frac{\pi}{2}$ there exist constants $C_{\varphi_{-},\varphi_{+}, k}$ and $\mu_0(\varphi_{-},\varphi_{+}, k)$ such that estimates (\ref{eq: eigenvalue approx}),(\ref{eq: int approx by Airy}) hold for all $ \varphi_0 \in [\varphi_{-},\varphi_{+}]$, $\mu \geq \mu_0(\varphi_{-},\varphi_{+}, k)$ and $C_{\varphi_{-},\varphi_{+},k}$ instead of $C$ . Indeed, this will be clear from the proof.
\end{rem}

Proposition \ref{prop: main approx} will follow from a repeated application of Lemma \ref{lem: perturbation}.

\begin{proof}[Proof of Proposition \ref{prop: main approx}]
We use the notation from the second half of Section \ref{subsec: perturbation estimates} but drop the label $R$ from the notation, i.e., $A$ denotes the Airy operator on $[0,\varphi_0\delta^{-1/3}]$, $u_k$ ($k \in \bbN$) its $L^2$-normalized eigenfunctions with Dirichlet boundary condition, and $\alpha_k$ ($k \in \bbN$) the eigenvalues. In contrast, $a_k$ are the eigenvalues of the Airy operator on $(0,\infty)$.

Throughout the proof, all implicit constants are understood to depend on $\varphi_0 \in (0,\frac{\pi}{2})$ (in fact, locally uniformly in $\varphi_0$ as in Remark \ref{rem: locally uniform constants}) and $k \in \bbN$, e.g.,  if we say "for all $\mu$ large enough" this is supposed to be understood as "for all $\mu$ large enough (depending only on $\varphi_0$ and $k$)". 

\smallskip\noindent
\textit{Claim 1. There exist $x^\ast >0$ and $C > 0$ such that for all $\mu$ large enough, we have
\[
	|u_k(x)| \leq Ce^{-x} 
	\quad \text{ for all } \quad 
	x \in [x^\ast,\varphi_0\delta^{-1/3}].
\]
In particular, for every polynomial $Q$ 
we have 
\[
	\int_0^{\varphi_0 \delta^{-1/3}}|Q(x)| u_k^2(x) \, dx \leq C_Q 
\]
for a finite constant $C_Q$ only depending on the degree of $Q$ and an upper bound for the absolute value of its coefficients (and $\varphi_0$ and $k$).}

\smallskip
The constants $x^\ast$ and $C$ depend on $k$, but in accordance with the convention mentioned before Claim 1, we do \emph{not} make this dependence explicit in the notation.

\smallskip\noindent
\textit{Proof of Claim 1.}
By definition $u_k^{\prime \prime}(x)=(x-\alpha_k)u_k(x)$ for all $x \in [0,\varphi_0\delta^{-1/3}]$. Thanks to Lemma \ref{lem: Airy on finite and infinite domain} we know $\alpha_k \leq a_k+1$ for all $\mu$ large enough. Thus, for all $x \geq x^\ast \coloneq a_k+2$ and $\mu$ large enough, the coefficient $x-\alpha_k$ is bounded from below by $1$. 

After potentially changing $u_k$ up to sign, we may without loss of generality assume that $u_k > 0$ close to $\varphi_0\delta^{-1/3}$. The Sturm comparison theorem shows for all $x \in [x^\ast,\varphi_0\delta^{-1/3}]$
\[
	u_k(x^\ast) \geq \cosh(x-x^\ast)u_k(x)+\sinh(x-x^\ast)(-u_k^\prime(x)) \geq \frac{1}{2}e^{x-x^\ast}u_k(x).
\]
Finally, $|u_k(x^\ast)| \leq C||u_k||_{H^1} \leq C||u_k||_{L^2} \leq C$ due to the Sobolev embedding, elliptic regularity, and the $L^2$-normalization. This establishes the exponential decay. The integral estimate follows from the exponential decay and the $L^2$-normalization.
\hfill$\blacksquare$

\smallskip\noindent
\textit{Claim 2. We have, for all $\mu$ large enough,
\[
	\tilde{\alpha}_k \leq \alpha_k + C \delta^{1/3}.
\]}

\smallskip\noindent
\textit{Proof of Claim 2.}
We will use the upper bound (\ref{eq: eigenvalue perturbation}). Towards doing so, we first observe that for all $u \in L^2\big(0,\varphi_0\delta^{-1/3}\big)$ we have
\[
  \Big|  ||u||_{\tilde{L}^2}^2-||u||_{L^2}^2 \Big| \leq \int_0^{\varphi_0 \delta^{-1/3}} \left|\frac{\cos^2(\varphi_0)}{\cos^2(\varphi_0 - \delta^{1/3}x)}-1 \right| u(x)^2 \, dx 
  \leq C \delta^{1/3}\int_0^{\varphi_0 \delta^{-1/3}} x u(x)^2 \, dx
\]
since $\cos^2(\varphi_0)/\cos^2(\cdot)$ is smooth, hence Lipschitz, on $[0,\varphi_0]$. Moreover, by linearity, Claim 1 holds for all $u \in {\rm span}\{u_1,\dots,u_k\}$ with $||u||_{L^2} \leq 1$. This implies 
\[
    \Big|  ||u||_{\tilde{L}^2}^2-||u||_{L^2}^2 \Big| \leq C \delta^{1/3}
\]
for all $u \in {\rm span}\{u_1,\dots,u_k\}$ with $||u||_{L^2} = 1$. Thus, (\ref{eq: perturbation norm distortion 1}) holds for some $\varepsilon_k=O(\delta^{1/3})$. Furthermore, for all $\mu$ large enough, $\alpha_k$ is uniformly bounded from above by Lemma \ref{lem: Airy on finite and infinite domain}. Therefore, it remains to show $||(\tilde{A}-A)u_i||_{L^2}=O(\delta^{1/3})$ for all $i=1,\dots,k$ in order to deduce Claim 2 from (\ref{eq: eigenvalue perturbation}). This will follow from elementary calculations (though slightly tedious) and Claim 1.

We first show that $\tilde A$ is indeed a perturbation of $A=-\frac{d^2}{dx^2}+x$. For this, we recall its definition \eqref{eq: def of perturbation operator} and write 
    \[
    \frac{\delta^{-1/3}}{2\tan(\varphi_0)}\left(1-\frac{\cos^2(\varphi_0)}{\cos^2(\varphi_0-\delta^{1/3}x)} \right) = x + \tilde e(x),
    \]
where $\tilde e(x)$ is defined by the equation above. Thus, $\tilde A = \frac{\cos^2(\varphi_0-\delta^{1/3}x)}{\cos^2(\varphi_0)}\left(-\frac{d^2}{dx^2}+ x + \tilde e(x) \right)$. 

Now, by Taylor,
\begin{equation*}
	 \frac{\cos^2(\varphi_0)}{\cos^2(\varphi_0-\theta)} = 1-2 \tan(\varphi_0)\theta+O\big(\theta^2\big)
	 \quad \text{for all } \theta \in [0,\varphi_0],
\end{equation*}
whence for all $x \in [0,\varphi_0\delta^{-1/3}]$
\begin{equation}\label{eq: expansion of coefficient}
	\frac{\delta^{-1/3}}{2\tan(\varphi_0)}\left(1-\frac{\cos^2(\varphi_0)}{\cos^2(\varphi_0-\delta^{1/3}x)} \right)=x+O\big(\delta^{1/3}x^2\big),
\end{equation}
so $\tilde e(x) = O\big(\delta^{1/3} x^2\big)$. 

Keeping in mind the definition \eqref{eq: def of perturbation operator} of $\tilde{A}$, we are now ready to compute the integral
\begin{multline*}
	||(\tilde{A}-A)u_i||_{L^2}^2=\int_0^{\varphi_0 \delta^{-1/3}} \left(1-\frac{\cos^2(\varphi_0-\delta^{1/3}x)}{\cos^2(\varphi_0)}\right)^2 \left(u_i^{\prime \prime}(x) \right)^2 \, dx \\
	+ \int_0^{\varphi_0 \delta^{-1/3}}\left(x-\frac{\cos^2(\varphi_0-\delta^{1/3}x)}{\cos^2(\varphi_0)} (x+\tilde e(x)) \right)^2u_i^2(x) \, dx.
\end{multline*}

For the first integral, we recall that $u_i^{\prime \prime}(x)=(x-\alpha_i)u_i(x)$ and note that $\left(1-
\frac{\cos^2(\varphi_0-\cdot)}{\cos^2(\varphi_0)}\right)$ is smooth, and hence Lipschitz, on $[0,\varphi_0]$. Thus 
\begin{align*}
	\int_0^{\varphi_0 \delta^{-1/3}} \left(1-\frac{\cos^2(\varphi_0-\delta^{1/3}x)}{\cos^2(\varphi_0)}\right)^2 \left(u_i^{\prime \prime}(x) \right)^2 \, dx 
	&\leq  
	C \delta^{2/3}\int_0^{\varphi_0 \delta^{-1/3}} x^2(x-\alpha_i)^2 u_i^2(x) \, dx \\
	&\leq  C \delta^{2/3},
\end{align*}
where for the second inequality we invoke Claim 1 with $Q(x)=x^2(x-\alpha_i)^2$. 

For the second integral, we write
\[
    x-\frac{\cos^2(\varphi_0-\delta^{1/3}x)}{\cos^2(\varphi_0)} (x+\tilde e(x)) = x \left(1 -\frac{\cos^2(\varphi_0-\delta^{1/3}x)}{\cos^2(\varphi_0)} \right)  -\frac{\cos^2(\varphi_0-\delta^{1/3}x)}{\cos^2(\varphi_0)}\tilde  e(x).
    \]
We use again the fact that $\left(1-
\frac{\cos^2(\varphi_0-\cdot)}{\cos^2(\varphi_0)}\right)$ is Lipschitz, so that the first term on the right-hand side is $O\big(\delta^{1/3}x^2\big)$. The second term is also $O\big(\delta^{1/3}x^2\big)$ because $\tilde{e}(x)=O\big(\delta^{1/3}x^2\big)$ and because the function in front of $\tilde e(x)$ is bounded. Therefore, the second integral is bounded by 
\begin{align*}
	 \int_0^{\varphi_0 \delta^{-1/3}}\left(x-\frac{\cos^2(\varphi_0-\delta^{1/3}x)}{\cos^2(\varphi_0)} (x+\tilde e(x)) \right)^2u_i^2(x) \, dx 
	 & \leq  C \delta^{2/3} \int_0^{\varphi_0 \delta^{-1/3}}x^4 u_i^2(x) \, dx \\
	 & \leq  C \delta^{2/3},
\end{align*}
where in the second inequality we again made use of Claim 1. This completes the proof of Claim 2.
\hfill$\blacksquare$

\smallskip\noindent
\textit{Claim 3. Claim 1 also holds for $\tilde{u}_k$.}

\smallskip\noindent
\textit{Proof of Claim 3.}
We write 
$\tilde{u}_k^{\prime \prime}(x)=\coeffutilde(x)\tilde{u}_k(x)$, where by (\ref{eq: coeff of rescaled eigenfct})
\[
    \coeffutilde(x)
    =\frac{\delta^{-1/3}}{2\tan(\varphi_0)}\left(1-\frac{\cos^2(\varphi_0)}{\cos^2(\varphi_0-\delta^{1/3}x)} \right)-\frac{\cos^2(\varphi_0)}{\cos^2(\varphi_0-\delta^{1/3}x)}\tilde{\alpha}_k.
\]
It follows from Claim 2 and Lemma \ref{lem: Airy on finite and infinite domain} that $\tilde{\alpha}_k \leq a_k+\frac{1}{2}$ for all $\mu$ large enough. Also $\cos^2(\varphi_0)/\cos^2(\cdot) \leq 1$ on $[0,\varphi_0]$. Moreover, note that the first summand in the above formula is increasing in $x$, and it is $x+O(\delta^{1/3}x^2)$ by (\ref{eq: expansion of coefficient}). Hence, by monotonicity, we have for all $x \geq x^\ast \coloneq a_k+2$ and $\mu$ large enough
\[
	\coeffutilde (x) \geq \coeffutilde(x^\ast) \geq x^\ast+O\big(\delta^{1/3}(x^\ast)^2\big)-\left(a_k+\frac{1}{2}\right) \geq 1.
\]
The rest of the proof of Claim 1 now carries over without any changes. This implicitly uses that the norms $||\cdot||_{\tilde{L}^2}$ and $||\cdot||_{L^2}$ are equivalent, so the usual $L^2$-norm of $\tilde{u}_k$ is bounded from above by a constant only depending on $\varphi_0$ and $k$.
\hfill$\blacksquare$

\smallskip\noindent
\textit{Claim 4. We have, for all $\mu$ large enough, 
\[
	|\tilde{\alpha}_k-\alpha_k| \leq C \delta^{1/3}.
\]}

\smallskip\noindent
\textit{Proof of Claim 4.} Using Claim 3 instead of Claim 1, the arguments from the proof of Claim 2 carry over to show that (\ref{eq: perturbation norm distortion 2}) holds with $\tilde{\varepsilon}_k=O(\delta^{1/3})$ and $||(\tilde{A}-A)\tilde{u}_i||_{L^2}=O(\delta^{1/3})$ for all $i=1,\dots,k$. Therefore, (\ref{eq: eigenvalue perturbation}) implies the desired estimate since $\tilde{\alpha}_k$ is uniformly bounded from above by Claim 2.
\hfill$\blacksquare$

\smallskip\noindent
\textit{Claim 5. After potentially changing $u_k$ up to sign, we have, for all $\mu$ large enough, 
\[
	||\tilde{u}_k-u_k||_{L^2} \leq C \delta^{1/3}.
\]}

\smallskip\noindent
\textit{Proof of Claim 5.} We know $||(\tilde{A}-A)u_i||_{L^2}=O(\delta^{1/3})$ and $||(\tilde{A}-A)\tilde{u}_i||_{L^2}=O(\delta^{1/3})$ for all $i=1,\dots,k$ from the proofs of Claim 2 and Claim 4. Moreover, by Lemma \ref{lem: Airy on finite and infinite domain}, all eigenvalues $\alpha_k$ of $A$ are simple and the distance between consecutive eigenvalues is uniformly bounded from below if $R=\varphi_0\delta^{-1/3}$ is large enough. Therefore, the desired estimate follows from (\ref{eq: eigenvector perturbation}).
\hfill$\blacksquare$

With Claims 1--5 at hand, we can now complete the proof of Proposition \ref{prop: main approx}.

Recall that $|\alpha_k-a_k|$ is exponentially small in $R=\varphi_0\delta^{-1/3}$ by Lemma \ref{lem: Airy on finite and infinite domain}. Together with Claim 4 this immediately yields Proposition \ref{prop: main approx}(i).

Similarly, $||u_k-v_k||_{L^2}$ is exponentially small in $R=\varphi_0\delta^{-1/3}$ by Lemma \ref{lem: Airy on finite and infinite domain}, and hence $||\tilde{u}_k-v_k||_{L^2}=O(\delta^{1/3})$ because of Claim 5. For any polynomial $Q$, we have
\[
    \int_0^{\varphi_0 \delta^{-1/3}}|Q(x)| \left|\tilde{u}_k^2(x)-v_k^2(x)\right| \, dx
    \leq \big(\left|\left|Q\tilde{u}_k \right|\right|_{L^2}
    +\left|\left|Qv_k \right|\right|_{L^2}\big)||\tilde{u}_k-v_k||_{L^2}
\]
due to the third binomial formula, the Cauchy-Schwarz inequality, and the triangle inequality. Note that, by the exact same proof, Claim 1 also holds for $v_k$.  Therefore the weighted integral estimate Proposition \ref{prop: main approx}(ii) follows from Claim 1 for $v_k$, Claim 3 and the unweighted estimate $||\tilde{u}_k-v_k||_{L^2}=O(\delta^{1/3})$.
\end{proof}

\subsection{Proof of Theorem \ref{Main Theorem} when \(n=2\)}\label{subsec: proof of main thm}

With Proposition \ref{prop: main approx}, we now complete the proof of Theorem \ref{Main Theorem} when $n=2$. We start by verifying the second assumption (\ref{eq: core of strategy 2}) in Lemma \ref{lem: strategy} for all $\mu$ large enough.

\begin{cor}\label{cor: integral expansion}
For every $\varphi_0 \in (0,\frac{\pi}{2})$ and $P \in C^2([0,\frac{\pi}{2}); \mathbb R)$ with $P'(\varphi_0)>0$ there exist constants $C > 0$ and $\mu_0 > 0$ with the following property. Let ${\rm I}$ denote the integral in (\ref{eq: core of strategy 2}), i.e.,
\[
	{\rm I} \coloneq  \int_0^{\varphi_0} P(\varphi)\left(h_2^2(\varphi)-h_1^2(\varphi)\right)\cos^{-2}(\varphi) \, d\varphi,
\]
where $h_1$ and $h_2$ are as in Lemma \ref{lem: strategy}. Then, for all $\mu \geq \mu_0$, we have
\[
	\left|\frac{\delta^{-1/3}}{P^\prime(\varphi_0)}{\rm I} + \frac{2}{3}(a_2-a_1) \right| \leq C \delta^{1/3},
\]
where $0 > -a_1 > -a_2 > \dots$ are the zeros of the Airy function $\Ai$. In particular, ${\rm I} < 0$ for all $\mu$ large enough.
\end{cor}

The constants $C$ and $\mu_0$ can be chosen to be locally uniform in $\varphi_0$ in the sense of Lemma \ref{rem: locally uniform constants}. Here $k=1,2$; so the constants are independent of $k$.

\begin{proof}
Due to the weighted $L^2$-normalization (\ref{eq: L^2 normalization h}) and the definition (\ref{eq: def of tilde(u)}) of $\tilde{u}_k$, we write
\begin{align*}
	{\rm I} \coloneq & \int_0^{\varphi_0} P(\varphi)\left(h_2^2(\varphi)-h_1^2(\varphi)\right)\cos^{-2}(\varphi) \, d\varphi \\
	=& \int_0^{\varphi_0} \big(P(\varphi)-P({\varphi_0})\big)\left(h_2^2(\varphi)-h_1^2(\varphi)\right)\cos^{-2}(\varphi) \, d\varphi 	\\
	=& \int_0^{{\varphi_0} \delta^{-1/3}} \big(P({\varphi_0}-\delta^{1/3}x)-P({\varphi_0})\big)\left(\frac{\tilde{u}_2^2(x)}{n_2^2}-\frac{\tilde{u}_1^2(x)}{n_1^2}\right)\cos^{-2}({\varphi_0}-\delta^{1/3}x) \, \delta^{1/3} dx .
\end{align*}
By Taylor $P({\varphi_0}-\delta^{1/3}x)-P({\varphi_0})=-P^{\prime}({\varphi_0})\delta^{1/3}x+O(\delta^{2/3}x^2)$ for all $x \in [0,{\varphi_0}\delta^{-1/3}]$. Thus,
\[
	\frac{\delta^{-1/3}}{P^{\prime}({\varphi_0})}{\rm I}=\int_0^{{\varphi_0} \delta^{-1/3}} \big(-x+O(\delta^{1/3}x^2)\big)\left(\frac{\delta^{1/3}}{n_2^2}\tilde{u}_2^2(x)-\frac{\delta^{1/3}}{n_1^2}\tilde{u}_1^2(x)\right)\cos^{-2}({\varphi_0}-\delta^{1/3}x) \, dx.
\]
From (\ref{eq: norm approx}) we know $\delta^{1/3}/n_k^2=\cos^2({\varphi_0})$. On the other hand, $\cos^{-2}(\cdot)$ is smooth, hence Lipschitz, on $[0,{\varphi_0}]$; so $\cos^{-2}({\varphi_0}-\delta^{1/3}x)=\cos^{-2}({\varphi_0})+O(\delta^{1/3}x)$. Finally, by Claim 3 in the proof of Proposition \ref{prop: main approx}, $\int |Q(x)|\tilde{u}_k^2(x) dx \leq C_Q < \infty$ for any polynomial $Q $. Hence
\[
	\frac{\delta^{-1/3}}{P^{\prime}({\varphi_0})}{\rm I}=\int_0^{{\varphi_0} \delta^{-1/3}}(-x)\left(\tilde{u}_2^2(x)-\tilde{u}_1^2(x)\right) \, dx + O(\delta^{1/3}).
\]
Proposition \ref{prop: main approx}(ii) and Lemma \ref{lem: model integral} yield
\begin{align*}
	\int_0^{{\varphi_0} \delta^{-1/3}}x\left(\tilde{u}_2^2(x)-\tilde{u}_1^2(x)\right) \, dx 
	=& \int_0^{{\varphi_0} \delta^{-1/3}}x\left(v_2^2(x)-v_1^2(x)\right) \, dx + O(\delta^{1/3}) \\
	=& \frac{2}{3}(a_2-a_1)+O(\delta^{1/3}).
\end{align*}
This completes the proof.
\end{proof}

Finally, we can present the proof of Theorem \ref{Main Theorem} when $n=2$.

\begin{proof}[Proof of Theorem \ref{Main Theorem} when $n=2$]
Choose $P$ as in \eqref{eq: def of potential}. For every $D \geq 0$ we denote by ${\varphi}_D \in (0,\frac{\pi}{2})$ the unique number such that the geodesic $[0,{\varphi}_D] \ni \varphi \mapsto (\sin(\varphi),\cos(\varphi)) \in \bbH^2$ has length $D$. 

Fix a diameter $D_0 > 0$. Set ${\varphi}_{\pm} \coloneq {\varphi}_{(1 \pm 1/2)D_0}$. Invoking Corollary \ref{cor: integral expansion} and Remark \ref{rem: locally uniform constants} we obtain a constant $\mu_0 > 0$ (only depending on $D_0$) with the following property: For every ${\varphi_0} \in [{\varphi}_{-},{\varphi}_{+}]$ and $\mu \geq \mu_0$ the integral ${\rm I}$ from Corollary \ref{cor: integral expansion} satisfies ${\rm I} < 0$, i.e., the second assumption (\ref{eq: core of strategy 2}) in Lemma \ref{lem: strategy} is satisfied. Similarly, after potentially increasing $\mu_0$, it follows from the asymptotic expansion (\ref{eq: eigenvalue asymptotic expansion}) that the first assumption (\ref{eq: strategy eigenvalue assumption}) in Lemma \ref{lem: strategy} is also satisfied for all ${\varphi_0} \in [{\varphi}_{-},{\varphi}_{+}]$ and $\mu \geq \mu_0$. Therefore, for all ${\varphi_0} \in [{\varphi}_{-},{\varphi}_{+}]$ and $\mu \geq \mu_0$, in $\Omega_{{\varphi},\mu}$ the constant potentials do \textit{not} minimise the spectral gap thanks to Lemma \ref{lem: strategy}. It remains to find one of these regions whose diameter equals $D_0$.

Note that for all $D > 0$ and $\mu > 0$
\begin{equation}\label{eq: diam estimates}
	D \leq {\rm diam}(\Omega_{{\varphi}_D,\mu}) \leq D + \cosh(D)\frac{\pi}{\sqrt{\mu}}.
\end{equation}
Indeed, this follows from the definition (\ref{eq: def of Omega}) of $\Omega_{{\varphi},\mu}$, the definition of ${\varphi}_D$, the triangle inequality, and hyperbolic geometry. 

After potentially increasing $\mu_0$ we can assume $D_0/2+\cosh(D_0/2)\pi/\sqrt{\mu_0} < D_0$. So, by (\ref{eq: diam estimates}), ${\rm diam}(\Omega_{{\varphi}_{-},\mu_0}) < D_0$ and ${\rm diam}(\Omega_{{\varphi}_{+},\mu_0}) > D_0$. Since ${\varphi} \mapsto {\rm diam}(\Omega_{{\varphi},\mu_0})$ is clearly continuous, we find some ${\varphi}_0 \in [{\varphi}_{-},{\varphi}_{+}]$ such that ${\rm diam}(\Omega_{{\varphi}_0,\mu_0})=D_0$. 
\end{proof}

\begin{rem}\label{rem: other gaps}\normalfont
Given $K \in \mathbb{N}, K \geq 2$, our estimates are uniform for all $k \leq K$. From the eigenvalue asymptotic expansion \eqref{eq: eigenvalue asymptotic expansion}, we get that $\lambda_k^{(\mu)} < \lambda_1^{(4 \mu)}$ for all $k \leq K$ for $\mu$ large enough. Combining these remarks with the fact that the proof of Lemma \ref{lem: model integral} works when $2$ and $1$ are replaced by $k$ and $k-1$ (or any $k$ and $j$), we obtain a similar result for all eigenvalue gaps up to $K$ if $\mu$ is large enough. In other words, for every $n\geq 2$, every $K\geq 2$ and every $D>0$, there exist a convex domain $\Omega \subseteq \mathbb H^n$ and a convex function $V$ on $\Omega$ such that $\text{diam}(\Omega) = D$ and $\lambda_k(\Omega;V) - \lambda_{k-1}(\Omega;V) < \lambda_k(\Omega) - \lambda_{k-1}(\Omega)$ for all $2\leq k \leq K$. 
\end{rem}

\section{Higher dimensions}\label{sec: high dim}

Until now we focused on the case $n=2$ for clarity of presentation. The goal of this section is to explain the adjustments needed to obtain Theorem \ref{Main Theorem} for general $n \geq 2$.

\subsection{Domain and change of variables}\label{subsec: set up - high dim}

Fix $n \geq 2$. We will now explain that then, in comparison to (\ref{eq: eq for h}), the relevant ODE eigenvalue problem to study is
\begin{equation}\label{eq: eq for h - high dim}
	\frac{\cos^{n-2}(\varphi)}{\cos^{n-2}(\varphi_0)}\frac{d}{d\varphi}\left(\frac{\cos^{n-2}(\varphi_0)}{\cos^{n-2}(\varphi)}h^\prime(\varphi)\right)=\big(\mu - \lambda \cos^{-2}(\varphi)\big)h(\varphi)
	\text{ for } \varphi \in (0,{\varphi_0})
\end{equation}
with Dirichlet boundary condition $h(0)=h({\varphi_0})=0$, where ${\varphi_0} \in (0,\frac{\pi}{2})$ will be a fixed parameter and we will take $\mu \to \infty$. The  weighted $L^2$-normalization will be
\begin{equation}\label{eq: L^2 normalization h - high dim}
	\int_0^{\varphi_0} \cos^{-n}(\varphi) h^2(\varphi) \, d\varphi=1.
\end{equation}
Accordingly, we will have to show that the integral
\begin{equation}\label{eq: core of strategy 2 - high dim}
    {\rm I} \coloneq \int_0^{\varphi_0}P(\varphi)\Big(h_2^2(\varphi)-h_1^2(\varphi) \big) \cos^{-n}(\varphi) \, d\varphi
\end{equation}
is negative (compare this with (\ref{eq: core of strategy 2}) in Lemma \ref{lem: strategy}).

Fix $H \subseteq \bbH^n$ a totally geodesic copy of $\bbH^{n-1}$. Fixing a unit normal vector field $\nu$ along $H$, we obtain a canonical diffeomorphism $\bbR \times H \stackrel{\cong}{\to} \bbH^n, (s,p) \mapsto \exp_p(s\nu(p))$, in which $s$ denotes the signed distance to $H$. In these coordinates, $g_{\bbH^n}=ds^2+\cosh^2(s)g_H$, where $g_H$ is the hyperbolic metric on $H$.

For any convex domain $\Omega_H \subseteq H$ and $\ell > 0$, we consider the convex domain
\begin{equation}\label{eq: def of Omega - high dim}
	\Omega \coloneq \Omega_{\Omega_H,\ell} \coloneq \Big\{ \exp_p(s\nu(p)) \in \bbH^n \, | \, p \in \Omega_H \text{ and } s \in (0,\ell)  \Big\}.
\end{equation}
We again take the convex potential $P \coloneq d_{\bbH^n}(\cdot,H)$; note that this is just $P(s,p)=s$ in the coordinates $\bbR \times H \cong \bbH^n$ from above.

In the coordinates $\bbR \times H \cong \bbH^n$, the Laplacian is given by 
\[
	\Delta u=\partial_{ss}^2(u)+(n-1)\tanh(s)\partial_s(u)+\frac{1}{\cosh^2(s)}\big(\Delta_H u(s,\cdot)\big),
\]
where $\Delta_H u(s,\cdot)$ denotes the $H$-Laplacian of the function $u(s,\cdot) \colon H \to \bbR, p \mapsto u(s,p)$. This again allows for a separation of variables ansatz. Namely, for $u(s,p)=h(s)f(p)$, solving $-\Delta u = \lambda u$ is equivalent to simultaneously solving
\[
	h^{\prime \prime}(s)+(n-1)\tanh(s)h^\prime(s)=\left(\frac{\tilde{\mu}}{\cosh^2(s)}-\lambda \right) h(s)
	\quad \text{and} \quad
	-\Delta_H f=\tilde{\mu}f,
\]
where $\tilde{\mu} \in \bbR$ is a constant. For $u$ to have Dirichlet boundary condition, $h$ and $f$ need to have Dirichlet boundary condition; so a non-zero solution to the above equations can only exist if $\tilde{\mu}=\mu_i$ for some $i \in \bbN$, where $0 < \mu_1 \leq \mu_2 \leq \dots $ are the eigenvalues of $-\Delta_H$ on $\Omega_H$ with Dirichlet boundary values. 

Assume from now on $\lambda_2^{(\mu_1)} < \lambda_1^{(\mu_2)}$, where $\lambda_k^{(\mu_i)}$ denotes the $k$-th eigenvalue of the equation for $h$ with $\tilde{\mu}=\mu_i$. Then, as before, $\lambda_k(\Omega)=\lambda_k^{(\mu_1)}$ and $u_k(s,p)=h_k^{(\mu_1)}(s)f_1(p)$ for $k=1,2$.

We can make choices for $\Omega_H$ such that $\mu_1(\Omega_H) \to \infty$, $\mu_2(\Omega_H)/\mu_1(\Omega_H) \to {\rm const}> 1$, and ${\rm diam}(\Omega_H) \to 0$; it then follows from (\ref{eq: eigenvalue asymptotic expansion}) that the assumption $\lambda_2^{(\mu_1)} < \lambda_1^{(\mu_2)}$ will be satisfied. For example, one can show (see \cite[Lemma 4.2 in the Appendix]{NSW22}) that, for all $k \in \bbN$, $\varepsilon^2 \lambda_k(B^{\bbH^m}({\rm pt},\varepsilon)) \to \lambda_k(B^{\bbR^m}({\rm pt},1))$ as $\varepsilon \to 0$; so $\Omega_H = B^H(p_0,\varepsilon)$ with $\varepsilon \to 0$ has the desired properties.

From now on we simply write $\mu=\mu_1(\Omega_H)$, and consider the above eigenvalue problem for $h$ with $\tilde{\mu}=\mu$. The desired equation (\ref{eq: eq for h - high dim}) follows from
\[
	h^{\prime \prime}(s)+(n-1)\tanh(s)h^\prime(s)=\left(\frac{\tilde{\mu}}{\cosh^2(s)}-\lambda \right) h(s)
	\text{ for } s \in (0,\ell)
\]
by a change of variables. Namely, we implicitly define a diffeomorphism $\varphi \colon \bbR \to (-\frac{\pi}{2},\frac{\pi}{2})$ by solving the ODE $\frac{d \varphi}{ds}(s)=\cos(\varphi(s))$ and $\varphi(0)=0$. Then $\tanh(s)=\sin(\varphi(s))$ for all $s \in \bbR$ since both functions solve the ODE $y^\prime(s)=1-y^2(s)$ and $y(0)=0$. Consequently, $\cos(\varphi)=(1-\sin^2(\varphi))^{1/2}=(1-\tanh(s))^{1/2}=1/\cosh(s)$. So the above equation becomes
\[
	\partial_{ss}^2h+(n-1)\sin(\varphi)\partial_s h =\big(\mu \cos^2(\varphi)-\lambda \big) h.
\]
Note $\partial_s h = (\partial_\varphi h) \cos(\varphi) $ and $\partial_{ss}^2 h = (\partial_{\varphi \varphi}^2 h)\cos^2(\varphi)+(\partial_\varphi h)(-\sin(\varphi))\cos(\varphi)$ since, by definition, $\partial_s \varphi = \cos(\varphi)$. Plugging this in, and dividing by $\cos^2(\varphi)$, we obtain
\[
	h^{\prime \prime}(\varphi)+(n-2)\tan(\varphi)h^\prime(\varphi)=\big(\mu - \lambda \cos^{-2}(\varphi) \big) h(\varphi)
	\text{ for }\varphi \in (0,{\varphi_0}),
\]
where ${\varphi_0} \coloneq \varphi(\ell) \in (0,\frac{\pi}{2})$. Dividing this by $\cos^{n-2}(\varphi)$ yields
\[
	\frac{d}{d\varphi}\left(\frac{1}{\cos^{n-2}(\varphi)}h^\prime(\varphi)\right)=\frac{\mu - \lambda \cos^{-2}(\varphi)}{\cos^{n-2}(\varphi)}h(\varphi),
\]
from which the desired equation (\ref{eq: eq for h - high dim}) readily follows. The form of the normalization (\ref{eq: L^2 normalization h - high dim}) comes from the fact that $d{\rm vol}_{\bbH^n}=\cosh^{n-1}(s) ds \wedge d{\rm vol}_H=\cos^{-n}(\varphi) d\varphi \wedge d{\rm vol}_H$ in the coordinates $\bbH^n \cong \bbR \times H \cong (-\frac{\pi}{2},\frac{\pi}{2}) \times H$ from above.

\subsection{Perturbation arguments}\label{subsec: Airy approx - high dim}

We now explain the necessary changes to extend the arguments from Section \ref{subsec: Airy approx} to the new eigenvalue problem (\ref{eq: eq for h - high dim}).

Fix ${\varphi_0} \in (0,\frac{\pi}{2})$ and take $\mu > 0$ large. Let $\lambda_1 \leq \lambda_2 \leq \dots$ be the eigenvalues and $h_k$ ($k \in \bbN$) the eigenfunctions of the eigenvalue problem (\ref{eq: eq for h - high dim}), i.e., 
\begin{equation*}\label{eq: eq for h - high dim 2}
	\frac{\cos^{n-2}(\varphi)}{\cos^{n-2}({\varphi_0})}\frac{d}{d\varphi}\left(\frac{\cos^{n-2}({\varphi_0})}{\cos^{n-2}(\varphi)}h_k^\prime(\varphi)\right)=\big(\mu - \lambda_k \cos^{-2}(\varphi)\big)h_k(\varphi)
	 \text{ for } \varphi \in (0,{\varphi_0})
\end{equation*}
with Dirichlet boundary condition $h_k(0)=h_k({\varphi_0})=0$ and the $L^2$-normalization (\ref{eq: L^2 normalization h - high dim}). We are again interested in the behaviour of $\lambda_k$ and $h_k$ as $\mu \to \infty$.

We define $\delta$, $\tilde{u}_k$, $\tilde{\alpha}_k$ by exactly the same formulas (\ref{eq: def of delta}), (\ref{eq: def of tilde(u)}), (\ref{eq: def of tilde(lambda)}) as in Section \ref{subsec: Airy approx}, while $n_k$ is chosen such that $\tilde{u}_k$ satisfies the weighted normalization
\begin{equation}\label{eq: L2 normalization tilde(u) - high dim}
    \int_0^{\varphi_0 \delta^{-1/3}}\frac{\cos^n(\varphi_0)}{\cos^{n}(\varphi_0-\delta^{1/3}x)}\tilde{u}_k(x)^2\,dx = 1.
\end{equation}
Explicitly,
\[
    \frac{\delta^{1/3}}{n_k^2}=\cos^n(\varphi_0).
\]
We will show below that Proposition \ref{prop: main approx} holds for all $n \geq 2$. Then the arguments from Section \ref{subsec: proof of main thm} carry over word by word, thus establishing Theorem \ref{Main Theorem} for all $n \geq 2$. The only differences will be that the constants $C,\mu_0$ are now also allowed to depend on $n$, and that one considers the integral (\ref{eq: core of strategy 2 - high dim}) instead of (\ref{eq: core of strategy 2}).

It remains to prove Proposition \ref{prop: main approx} for $n \geq 2$. To achieve this, we need to adapt the definition of the operator $\tilde{A}$ that is close the Airy operator and that satisfies, for all $k \in \bbN$,
\[
	\tilde{A}\tilde{u}_k=\tilde{\alpha}_k\tilde{u}_k.
\]
Using the definition (\ref{eq: def of tilde(u)}) of $\tilde{u}_k$ and the equation \eqref{eq: eq for h - high dim} for $h_k$, we observe that
\begin{equation}\label{eq: eq for tilde(u)}
	\frac{\cos^{n-2}({\varphi_0}-\delta^{1/3}x)}{\cos^{n-2}({\varphi_0})}\frac{d}{dx}\left(\frac{\cos^{n-2}({\varphi_0})}{\cos^{n-2}({\varphi_0}-\delta^{1/3}x)}\tilde{u}_k^\prime(x)\right)=
    \coeffutilde
    (x)\tilde{u}_k(x),
\end{equation}
where by the exact same calculation as in (\ref{eq: coeff of rescaled eigenfct})
\begin{equation}\label{eq: coeff of rescaled eigenfct - high dim}
    \coeffutilde(x) 
	= \frac{\delta^{-1/3}}{2\tan({\varphi_0})}\left(1-\frac{\cos^2({\varphi_0})}{\cos^2({\varphi_0}-\delta^{1/3}x)} \right)-\frac{\cos^2({\varphi_0})}{\cos^2({\varphi_0}-\delta^{1/3}x)}\tilde{\alpha}_k.
\end{equation}
Therefore,
\begin{equation}\label{eq: def of perturbation operator - high dim}
\begin{split}
	\tilde{A} \coloneq & -\frac{\cos^{n}({\varphi_0}-\delta^{1/3}x)}{\cos^{n}({\varphi_0})}\frac{d}{dx}\left(\frac{\cos^{n-2}({\varphi_0})}{\cos^{n-2}({\varphi_0}-\delta^{1/3}x)}\frac{d}{dx}\right) \\
	& +\frac{\cos^2({\varphi_0}-\delta^{1/3}x)}{\cos^2({\varphi_0})}\frac{\delta^{-1/3}}{2\tan({\varphi_0})}\left(1-\frac{\cos^2({\varphi_0})}{\cos^2({\varphi_0}-\delta^{1/3}x)} \right)
\end{split}
\end{equation}
satisfies $\tilde{A}\tilde{u}_k=\tilde{\alpha}_k\tilde{u}_k$ for all $k \in \bbN$.

In the same way as previously, $\tilde{A}$ is an unbounded self-adjoint operator on the weighted $L^2$-space 
\[
    \tilde{L}^2(0,\varphi_0 \delta^{-1/3}) \coloneq \left(L^2(0,\varphi_0\delta^{-1/3}), \langle \cdot , \cdot \rangle_{\tilde{L}^2} \right)
\]
with the weighted $L^2$-inner product $\langle \cdot , \cdot \rangle_{\tilde{L}^2}$ given by 
\begin{equation}\label{eq: weighted L^2-norm}
    \langle f , g \rangle_{\tilde{L}^2} \coloneq \int_0^{\varphi_0 \delta^{-1/3}} \frac{\cos^n(\varphi_0)}{\cos^n(\varphi_0 - \delta^{1/3}x)} f(x)g(x) \, dx.
\end{equation}
The $\tilde{u}_k$ ($k \in \bbN$) then form a complete orthonormal basis for $\tilde{L}^2(0,\varphi_0 \delta^{-1/3})$ due to the weighted $L^2$ normalization (\ref{eq: L2 normalization tilde(u) - high dim}).

\begin{proof}[Proof of Proposition \ref{prop: main approx} when $n \geq 2$]
Throughout the proof we will use the same conventions and notations as in the proof of Proposition \ref{prop: main approx} in Section \ref{subsec: Airy approx}, the only difference being that the implicit constants are now also allowed to depend on $n$. 

\smallskip\noindent
\textit{Claim 1. There exist $x^\ast >0$ and $C > 0$ such that for all $\mu$ large enough we have
\[
	|u_k|(x),|u_k^\prime|(x) \leq Ce^{-x} 
	\quad \text{ for all } \quad 
	x \in [x^\ast,{\varphi_0}\delta^{-1/3}].
\]
In particular, for every polynomial $Q$ 
we have 
\[
	\int_0^{{\varphi_0} \delta^{-1/3}}|Q|(x) u_k^2(x) \, dx \leq C_Q 
	\quad \text{and}\quad
	\int_0^{{\varphi_0} \delta^{-1/3}}|Q|(x) \big(u_k^\prime(x)\big)^2 \, dx \leq C_Q 
\]
for a finite constant $C_Q$ only depending on $Q$ (and ${\varphi_0}$, $k$, and $n$).}

\smallskip\noindent
\textit{Proof of Claim 1.}
This follows directly from the proof of Claim 1 in Section \ref{subsec: Airy approx}. 
\hfill$\blacksquare$

\smallskip\noindent
\textit{Claim 2. We have, for all $\mu$ large enough,
\[
	\tilde{\alpha}_k \leq \alpha_k + C \delta^{1/3}.
\]}

\smallskip\noindent
\textit{Proof of Claim 2.}
We will only show that $\tilde{A}$ is a perturbation of $A$. Then Claim 2 can be deduced from (\ref{eq: eigenvalue perturbation}) exactly as in Section \ref{subsec: Airy approx}.

From the definition (\ref{eq: def of perturbation operator - high dim}) of $\tilde{A}$ we have
\begin{align*}
	\tilde{A}=& -\frac{\cos^{2}({\varphi_0}-\delta^{1/3}x)}{\cos^{2}({\varphi_0})}\frac{d^2}{dx^2} \\
	&-\frac{\cos^{n}({\varphi_0}-\delta^{1/3}x)}{\cos^{n}({\varphi_0})}\frac{d}{dx}\left(\frac{\cos^{n-2}({\varphi_0})}{\cos^{n-2}({\varphi_0}-\delta^{1/3}x)}\right)\frac{d}{dx} \\
	&+\frac{\cos^2({\varphi_0}-\delta^{1/3}x)}{\cos^2({\varphi_0})}\frac{\delta^{-1/3}}{2\tan({\varphi_0})}\left(1-\frac{\cos^2({\varphi_0})}{\cos^2({\varphi_0}-\delta^{1/3}x)} \right).
\end{align*}
Note $\cos^2({\varphi_0}-\delta^{1/3}x)/\cos^2({\varphi_0})=1+O(\delta^{1/3}x)$ since smooth functions are Lipschitz. Clearly, the coefficient of $\frac{d}{dx}$ is $O(\delta^{1/3})$ because of the chain rule. Thus, by (\ref{eq: expansion of coefficient}), 
\[
	\tilde{A}-A=O\big(\delta^{1/3}x \big)\frac{d^2}{dx^2}+
	O\big(\delta^{1/3}\big)\frac{d}{dx}+
	O\big(\delta^{1/3}x^2\big){\rm id}.
\]
The rest follows from Claim 1 exactly as in Section \ref{subsec: Airy approx}.
\hfill$\blacksquare$

\smallskip\noindent
\textit{Claim 3. There exists $x^\ast >0$ and $C > 0$ such that for all $\mu$ large enough we have
\[
	|\tilde{u}_k|(x) \leq Ce^{-x} 
	\quad \text{ for all } \quad 
	x \in [x^\ast,{\varphi_0}\delta^{-1/3}]
\]
and
\[
    \int_x^{\varphi_0 \delta^{-1/3}}(\tilde{u}_k^\prime)^2(\zeta) \, d\zeta \leq Ce^{-x} 
	\quad \text{ for all } \quad 
	x \in [x^\ast,{\varphi_0}\delta^{-1/3}].
\]
In particular, for every polynomial $Q$ 
we have 
\[
	\int_0^{{\varphi_0} \delta^{-1/3}}|Q|(x) \tilde{u}_k^2(x) \, dx \leq C_Q 
	\quad \text{and}\quad
	\int_0^{{\varphi_0} \delta^{-1/3}}|Q|(x) \big(\tilde{u}_k^\prime(x)\big)^2 \, dx \leq C_Q.
\]}

With a tiny bit work one can also show that $|\tilde{u}_k^\prime(x)| \leq Ce^{-x}$ for all $x \geq x^\ast$. Since this integral estimate is sufficient for our purposes, we refrain from doing so.

\smallskip\noindent
\textit{Proof of Claim 3.}
Exactly as in Section \ref{subsec: Airy approx} we deduce from Claim 2 that there exists $x^\ast > 0$ such that $\coeffutilde(x)\geq 4$ for all $x \geq x^\ast$ and $\mu$ large enough, where $\coeffutilde$ is given by (\ref{eq: coeff of rescaled eigenfct - high dim}). From (\ref{eq: eq for tilde(u)}) we recall
\[
	\tilde{u}_k^{\prime \prime}(x)+\frac{\cos^{n-2}({\varphi_0}-\delta^{1/3}x)}{\cos^{n-2}({\varphi_0})}\frac{d}{dx}\left(\frac{\cos^{n-2}({\varphi_0})}{\cos^{n-2}({\varphi_0}-\delta^{1/3}x)}\right)\tilde{u}_k^\prime(x)=
    \coeffutilde \tilde{u}_k(x).
\]
Thus, by the chain rule,
\[
	\tilde{u}_k^{\prime \prime}(x) = 
    \coeffutilde (x) \tilde{u}_k(x) + \delta^{1/3}\tilde{b}(x)\tilde{u}_k^\prime(x),
\]
for some continuous function $\tilde{b} : [0,{\varphi_0}\delta^{-1/3}] \to \bbR_{\geq 0}$ with $||\tilde{b}||_{C^0}=O(1)$.
Moreover, after potentially changing the sign of $\tilde{u}_k$, we may assume $\tilde{u}_k^\prime({\varphi_0}\delta^{-1/3}) < 0$. 

Observe that for all $x \geq x^\ast$ and all $\delta > 0$ small enough we have
\begin{align}
\notag
   \left( \frac{1}{2}\tilde{u}_k^2 \right)^{\prime \prime}
   = \left(\tilde{u}_k \tilde{u}_k^\prime \right)^\prime
   = \left(\tilde{u}_k^\prime\right)^2  + \tilde{u}_k\tilde{u}_k^{\prime \prime} 
   &=\left(\tilde{u}_k^\prime\right)^2  + \tilde{u}_k \left(\tilde{c}_k\tilde{u}_k+\delta^{1/3}\tilde{b} \tilde{u}_k^\prime \right) \\
   &\geq  \frac{1}{2}\left(\tilde{u}_k^\prime\right)^2 + \frac{\tilde{c}_k}{2}\tilde{u}_k^2 
   \geq  \sqrt{\tilde{c}_k}\left|\tilde{u}_k \tilde{u}_k^\prime\right|
   \geq 2 \left|\tilde{u}_k \tilde{u}_k^\prime \right|. \label{eq:tildeu_ksecond}
\end{align}
In particular, $\left( \frac{1}{2}\tilde{u}_k^2 \right)^{\prime \prime} \geq \frac{\tilde{c}_k}{2}\tilde{u}_k^2 \geq 4\left( \frac{1}{2}\tilde{u}_k^2 \right)$ for all $x \geq x^\ast$. Thus, by the Sturm comparison theorem, for all $x \in [x^\ast,\delta^{-1/3}\varphi_0]$ we get
\[
   \tilde{u}_k^2(x^\ast) \geq \cosh(2(x-x^\ast))\tilde{u}_k^2(x) \geq \frac{1}{2}e^{2(x-x^\ast)}\tilde{u}_k^2(x).
\]
Keeping in mind that $|\tilde{u}_k|(x^\ast) \leq C$ due to the Sobolev embedding, elliptic regularity, and the $L^2$-normalization (\ref{eq: L^2 normalization h - high dim}), this establishes the exponential decay of $\tilde{u}_k$. 

It remains to prove the estimate for $\tilde{u}_k^\prime$. For this we write $\hat{u}_k(t) \coloneq \tilde{u}_k(\delta^{-1/3}\varphi_0-t)$ and $t^\ast \coloneq \delta^{-1/3}\varphi_0-x^\ast$. From \eqref{eq:tildeu_ksecond} we then have, $(\hat{u}_k \hat{u}_k^\prime)^\prime(t) \geq 2 |\hat{u}_k \hat{u}_k^\prime|(t)$ for all $t \in [0,t^\ast]$. From this one readily deduces, for all $t \in [0,t^\ast]$, 
\[
    \hat{u}_k(t^\ast) \hat{u}_k^\prime(t^\ast) \geq e^{2(t^\ast-t)}  \hat{u}_k(t) \hat{u}_k^\prime(t),
\]
and thus $\hat{u}_k(t) \hat{u}_k^\prime(t) \leq Ce^{-2(t^\ast-t)}$ due to the Sobolev embedding, elliptic regularity, and the $L^2$-normalization (\ref{eq: L^2 normalization h - high dim}).
Moreover, from \eqref{eq:tildeu_ksecond} we also have $\left(\hat{u}_k \hat{u}_k^\prime \right)^\prime \geq \frac{1}{2}\left(\hat{u}_k^\prime\right)^2$ for all $t \in [0,t^\ast]$. Integrating this from $0$ to $t$ implies
\[
    \frac{1}{2}\int_0^{t} \left(\hat{u}_k^\prime\right)^2(\tau) \, d\tau \leq \hat{u}_k(t) \hat{u}_k^\prime(t) \leq C e^{-2(t^\ast-t)}.
\]
Since $t^\ast-t=x-x^\ast$ and $x^\ast$ is uniformly bounded, this completes the proof of Claim 3.
\hfill$\blacksquare$

The rest of the arguments from the proof of Proposition \ref{prop: main approx} in Section \ref{subsec: Airy approx} now carry over without any difficulties.
\end{proof}


\end{document}